\documentclass[11pt]{amsart}
\usepackage[top=1 in,bottom=1in, left=1 in, right=1 in]{geometry}
\usepackage[maxbibnames=99, style=alphabetic, maxalphanames=99, minalphanames=6, 
sorting=nyt, giveninits=true, useprefix=true]{biblatex}
\AtEveryBibitem{
    \clearfield{issn}
  \clearfield{isbn}
  \clearfield{doi}
  \clearlist{location}
  \clearfield{month}
  \clearfield{day}
  \clearfield{series}
  \clearlist{language}
  \clearlist{translator}
  \renewbibmacro{in:}{}
}
\DeclareLabelalphaNameTemplate{
  \namepart[use=true, pre=true, strwidth=1, compound=true]{prefix}
  \namepart[compound=true]{family}
}

\usepackage{amssymb, mathtools, thmtools}
\usepackage{mathdots}
\usepackage{graphicx}
\usepackage{microtype}
\usepackage{accents}
\usepackage{enumitem}
\usepackage{indentfirst}
\usepackage{xcolor}
\usepackage{tikz-cd}
\usepackage{csquotes}
\usepackage[colorlinks=true, linkcolor=blue]{hyperref} 
\usepackage[noabbrev]{cleveref}
\usepackage{bookmark}

\makeatletter
\ifdefined\newcounteralias 
  \let\fixthm@newcounteralias\newcounteralias
  \renewcommand*\newcounteralias[2]{%
    \def\fixthm@a{#1}\def\fixthm@b{#2}%
    \ifx\fixthm@a\fixthm@b\else\fixthm@newcounteralias{#1}{#2}\fi}
\fi
\makeatother

\theoremstyle{definition}
\declaretheorem[name=Definition,refname={definition,definitions},Refname={Definition,Definitions},numberwithin=section]{definition}

\theoremstyle{plain}
\declaretheorem[name=Lemma,refname={lemma,lemmas},Refname={Lemma,Lemmas},sibling=definition]{lemma}

\declaretheorem[name=Theorem,refname={theorem,theorems},Refname={Theorem,Theorems},sibling=definition]{thm}
\declaretheorem[name=Corollary,refname={corollary,corollaries},Refname={Corollary,Corollaries},sibling=definition]{corollary}
\theoremstyle{remark}

\newcommand{\GL}{\mathrm{GL}}
\newcommand{\numberset}{\mathbb}

\newcommand{\C}{\numberset{C}}

\newcommand{\g}{\mathfrak{g}}
\newcommand{\git}{/\!\!/}

\newcommand{\overbar}[1]{\mkern 1.5mu\overline{\mkern-1.5mu#1\mkern-1.5mu}\mkern 1.5mu}

\DeclareMathOperator{\id}{Id}

\newcommand{\hh}{\mathfrak{h}}
\newcommand{\Wh}{\mathrm{Wh}}
\DeclareMathOperator{\Hom}{Hom}
\DeclareMathOperator{\spec}{Spec}

\title{Symplectic singularities in Relative Langlands duality for products of general linear groups}
\date{}
\author{Bruno da Silveira Dias}
\author{Victor Ginzburg}

\begin{document}
	
    \begin{abstract}
        It was shown by Braverman-Finkelberg-Nakajima that the regular perverse sheaf on the affine Grassmannian of type A is isomorphic to the ring object associated to a Gelfand-Tsetlin quiver. We use this result to identify, for any finite dimensional representation $\mathbf{N}$ of a product of general linear groups, the relative Langlands dual $(T^*\mathbf{N})^\vee$ of $T^*\mathbf{N}$ with the Coulomb branch of a suitable representation. By a result of Bellamy, this implies that $(T^*\mathbf{N})^\vee$ has symplectic singularities.
    \end{abstract}
    \maketitle
    
    \section{Introduction}

    Let $G$ be a complex connected reductive group and $G^\vee$ the Langlands dual group. In \cite{benzvi2024relativelanglandsduality}, Ben-Zvi, Sakellaridis, and Venkatesh proposed a relative Langlands duality conjecture relating certain affine Hamiltonian $G$-varieties to affine Hamiltonian $G^\vee$-varieties. In the case where the Hamiltonian $G$-variety has the form $T^*\mathbf{N}$, for a finite dimensional $G$-representation $\mathbf{N}$, the dual variety $(T^*\mathbf{N})^\vee$ is obtained from the Satake transform of a ring object $\mathcal{A}_{G,\mathbf{N}}$ in the derived Satake category; see \cite{BFN-Ring}. We show that when $G$ is a product of general linear groups the variety $(T^*\mathbf{N})^\vee$ is isomorphic to the Coulomb branch associated to a representation of a larger group that contains $G$ as a subgroup. This implies, thanks to a result of Bellamy \cite{Bellamy2023}, that $(T^*\mathbf{N})^\vee$ has symplectic singularities. We also establish similar results for the affinization of the \emph{transfer variety} $\Wh_{H \to G}$ that plays an important role in functoriality for geometric Satake; see \cite[Section 2]{ginzburg2025pointwisepurityderivedsatake}.

    \subsection{Statement of Results}
    For each $n \geq 1$, let
    \[
        \mathcal{G}_n := \prod_{i=1}^n \GL_i,\quad  \mathcal{V}_n := \bigoplus_{i=1}^{n-1} \Hom(\C^i,\C^{i+1}).
    \]
    These are the group and representation space associated to the Gelfand-Tsetlin quiver
    \[
        \C^1 \to \C^2 \to \dots \to \C^n.
    \]
    We let $p_n\colon \mathcal{G}_n \to \GL_n$ denote the projection onto the last factor.

    Let $G := \prod_{j=1}^s\GL_{n_j}$ be a product of general linear groups, for some $n_j\geq 1$. Set $\mathcal{G} := \prod_{j=1}^s\mathcal{G}_{n_j}$, $\mathcal{V}:=\bigoplus_{j=1}^s\mathcal{V}_{n_j}$, and let $p\colon \mathcal{G} \to G$ be the product of the projections $p_{n_j}$. Our first result identifies the relative Langlands dual of $T^*\mathbf{N}$, for any representation $\mathbf{N}$ of $G$, with the Coulomb branch of a suitable representation of $\mathcal{G}$.

    \begin{thm}\label{thm-1}
        For any finite-dimensional representation $\mathbf{N}$ of $G$, there is a  $G$-equivariant isomorphism of Poisson schemes
        \[
            (T^*\mathbf{N})^\vee \cong \mathcal{M}_C(\mathcal{G}, \mathcal{V} \oplus p^*\mathbf{N}).
        \]
    \end{thm}

    Our second result uses a similar argument to provide a Coulomb branch description of the affine varieties $\overbar{\mathrm{Wh}_{H \to G}}$ that occur in the functoriality for derived Satake established in \cite[Section 2]{ginzburg2025pointwisepurityderivedsatake}. Before stating our result precisely, we recall their construction. 
    
    Given a reductive group $G$, denote its Langlands dual group by $G^\vee$. Fix a maximal unipotent subgroup $U^\vee \subseteq G^\vee$ and a non-degenerate character $\psi\colon \mathfrak{u}^\vee \to \C$. Let 
    \[
        \mathrm{Wh}_{G^\vee} := T^*_{\psi}(G^\vee/U^\vee) = G^\vee \times_{U^\vee} (\psi + (\mathfrak{u}^\vee)^\perp)
    \]
    denote the Whittaker-twisted cotangent bundle of $G^\vee/U^\vee$. It is a smooth affine symplectic variety equipped with commuting Hamiltonian actions of $G^\vee$ and of the universal centralizer $J_{G^\vee}$. As explained in \cite{ginzburg2025pointwisepurityderivedsatake}, a morphism $\rho\colon H \to G$ induces a map from $\mathfrak{c}_{H^\vee} \cong \hh\git H \to \mathfrak{c}_{G^\vee} \cong \g\git G$ and a morphism from $\rho^*J_{G^\vee}$ to $J_{H^\vee}$ of group schemes over $\mathfrak{c}_{H^\vee}$. Taking $\rho\colon H \to H \times G$, the graph embedding of a morphism $\phi\colon H \to G$, one defines the quasi-affine scheme
    \[
        \mathrm{Wh}_{H \to G} := J_{H^\vee} \times_{\mathfrak{c}_{H^\vee}}^{\rho^*J_{H^\vee \times G^\vee}} \rho^*\mathrm{Wh}_{H^\vee \times G^\vee} = \mathrm{Wh}_{H^\vee} \times_{\mathfrak{c}_{H^\vee}}^{\phi^*J_{G^\vee}} \phi^*\mathrm{Wh}_{G^\vee} .
    \]
    It is shown in \cite[Theorem 2.1.2]{ginzburg2025pointwisepurityderivedsatake} that the ring of functions on this scheme is a finitely generated integrally closed domain whose spectrum, which we denote by $\overbar{\mathrm{Wh}_{H \to G}}$, is an affine Poisson variety with a Hamiltonian $H^\vee \times G^\vee$-action.
     
    Now assume 
    \[
        H := \prod_{i=1}^r\GL_{m_i},\quad G := \prod_{j=1}^{s}\GL_{n_j}
    \]
    and let $\phi\colon H \to G$ be a homomorphism. Put $\mathcal{H}:=\prod_{i=1}^r\mathcal{G}_{m_i}$, $\mathcal{G} :=\prod_{j=1}^s\mathcal{G}_{n_j}$, $\mathcal{K}:=\prod_{j=1}^s\mathcal{G}_{n_j-1}$, and let $q\colon \mathcal{H} \to H$, resp.\ $p\colon \mathcal{G} \to G$, be the product of the projections $p_{m_i}$, resp.\ $p_{n_j}$. Identifying $\mathcal{K} \times G = \mathcal{G}$, consider the embedding $\iota_\phi\colon \mathcal{H} \times \mathcal{K} \hookrightarrow \mathcal{H} \times \mathcal{G}$ given by
    \begin{equation}\label{equation-homomorphism1}
        (h,k) \mapsto (h,k,\phi(q(h))).
    \end{equation}
    Write $\mathcal{U}:=\bigoplus_{i=1}^r\mathcal{V}_{m_i}$ and $\mathcal{V}:=\bigoplus_{j=1}^{s}\mathcal{V}_{n_j}$.

    \begin{thm}\label{thm-2}
        There is an isomorphism of affine schemes
        \[
            \mathcal{M}_C(\mathcal{H} \times \mathcal{K}, \iota_\phi^*(\mathcal{U} \oplus \mathcal{V})) \cong \overbar{\mathrm{Wh}_{H \to G}}.
        \]
    \end{thm}

    From a theorem of Bellamy \cite{Bellamy2023}, we obtain the following:

    \begin{corollary}
        For any finite-dimensional representation $\mathbf{N}$ of $G = \prod_{j=1}^s\GL_{n_j}$, the relative Langlands dual of $T^*\mathbf{N}$ has symplectic singularities. 
        
        Moreover, if $H = \prod_{i=1}^r\GL_{m_i}$ then for any algebraic group homomorphism $\phi\colon H \to G$, the affine variety $\overbar{\mathrm{Wh}_{H \to G}}$ has symplectic singularities.
    \end{corollary}

    \subsection*{Acknowledgments} The proofs given below are the product of useful conversations with ChatGPT-6, which was also used for proofreading. All writing was done by the authors.

    \section{Proofs}
    
    \subsection{Preliminaries} Given $G$ a complex reductive group and $\mathbf{N}$ a finite-dimensional representation of $G$, we denote by $\mathcal{A}_{G,\mathbf{N}}$ the corresponding ring object in the ind-completed derived Satake category $D_{G(O)}(\mathrm{Gr}_G)$, as defined in \cite{BFN2018,BFN-Ring}. The Coulomb branch associated to $(G,\mathbf{N})$ is the affine scheme
    \[
        \mathcal{M}_C(G,\mathbf{N}) := \spec(H_{G(O)}^\bullet(\mathrm{Gr}_G,\mathcal{A}_{G,\mathbf{N}})).
    \]

    Let $\mathcal{R}_G \in D_{G(O)}(\mathrm{Gr}_G)$ denote the regular sheaf, i.e.\ the object corresponding to the regular representation $\C[G^\vee]$ under geometric Satake equivalence. According to \cite[Lemma 5.13(c)]{BFN-Ring} the derived Satake equivalence $\Phi$ may be defined as the functor $\mathcal{F}\mapsto R\Gamma_{G(O)}(\mathrm{Gr}_G, \mathcal{R}_G \otimes^! \mathcal{F})$, cf.\ also \cite[Section 5]{ginzburg2025pointwisepurityderivedsatake}. Here $\mathcal{F}_1\otimes^!\mathcal{F}_2:= \Delta^!(\mathcal{F}_1\boxtimes\mathcal{F}_2)$, where $\Delta$ is the diagonal embedding. By \cite[Proposition 5.15(b)]{BFN-Ring}, for a ring object $\mathcal{F}$, this identification is compatible with the the multiplications. For a BFN ring object $\mathcal{A}_{G,\mathbf{N}}$,
    \[
        H^\bullet(\Phi(\mathcal{A}_{G,\mathbf{N}})) \cong H_{G(O)}^\bullet(\mathrm{Gr}_G, \mathcal{R}_G \otimes^! \mathcal{A}_{G,\mathbf{N}})
    \]
    is a commutative algebra over $\mathrm{Sym}(\g^\vee)$. Commutativity follows from \cite[Proposition 2.5]{BFN-Ring}, \cite[Lemmas 2.2 and 2.5-2.6]{BFN2018}, and \cite{MirkovicVilonen2007}. In this setting, by definition, the relative Langlands dual of $T^*\mathbf{N}$ is
    \[
        (T^*\mathbf{N})^\vee:= \spec H_{G(O)}^\bullet(\mathrm{Gr}_G, \mathcal{R}_G \otimes^! \mathcal{A}_{G,\mathbf{N}}).
    \]

    The following result is a consequence of \cite[Theorem 2.11]{BFN-Ring}; see \cite[proof of Lemma 3.5]{gannon2025functorialitycoulombbranches}. Another proof can be found in \cite{dias2026rankinselbergdualitygluing}.

    \begin{lemma}\label{lemma-gt}
        There is an isomorphism of ring objects $(p_n)_*\mathcal{A}_{\mathcal{G}_n, \mathcal{V}_n} \cong \mathcal{R}_{\GL_n}$. 
    \end{lemma}

    We will also use the following general facts. The first is proved in \cite[3(vii)]{BFN2018}, the second is stated in \cite[5(xi)]{BFN-Ring}, and the third is \cite[Proposition 3.3]{GannonWilliams2026}:
    \begin{enumerate}[label=\roman*)]
        \item If $(G,\mathbf{N}) = (G_1 \times G_2, \mathbf{N}_1 \oplus \mathbf{N}_2)$, where $\mathbf{N}_i$ is a representation of $G_i$, then
        \begin{equation}\label{equation-fact1}
            \mathcal{A}_{G,\mathbf{N}} = \mathcal{A}_{G_1,\mathbf{N}_1} \boxtimes \mathcal{A}_{G_2,\mathbf{N}_2}.
        \end{equation}
        
        \item If $\rho\colon H \hookrightarrow G$ is a closed embedding of reductive groups and $\mathbf{N}$ is a representation of $G$, then
        \begin{equation}\label{equation-fact2}
            \rho^\dagger\mathcal{A}_{G,\mathbf{N}} \cong \mathcal{A}_{H,\rho^*\mathbf{N}}.
        \end{equation}

        \item Consider a cartesian square of connected reductive groups
        \[
            \begin{tikzcd}
                \widetilde{H} \arrow[r, "\widetilde{\rho}"]\arrow[d, "q_H"'] & \widetilde{G} \arrow[d, "q_G"] \\
                H \arrow[r, "\rho"] & G
            \end{tikzcd}
        \]
        and assume $q_G$ and $q_H$ are quotient maps. If $\mathcal{A} \in D_{\widetilde{G}(O)}(\mathrm{Gr}_{\widetilde{G}})$ is a ring object then there is a canonical isomorphism of ring objects
        \begin{equation}\label{equation-square}
            (q_H)_*\widetilde{\rho}^\dagger\mathcal{A} \cong \rho^\dagger(q_G)_*\mathcal{A}
        \end{equation}
        where $(-)^\dagger$ denotes restriction of equivariance followed by $!$-pullback.
    \end{enumerate}
    
    \subsection{Proof of \texorpdfstring{\Cref{thm-1}}{Theorem 1.1}}

    Consider the diagram
    \[
        \begin{tikzcd}
            \mathcal{G} \arrow[r, "\iota"]\arrow[d, "p"'] & \mathcal{G} \times G \arrow[d, "p \times \id"] \\
            G \arrow[r, "\Delta"] & G \times G
        \end{tikzcd}
    \]
    where $\iota$ is the graph embedding $\iota(g)=(g,p(g))$ and $\Delta$ is the diagonal map. Apply \eqref{equation-square} to the square above and the ring object $\mathcal{A}_{\mathcal{G} \times G,\mathcal{V} \oplus \mathbf{N}}$. By \eqref{equation-fact2}, $\iota^\dagger\mathcal{A}_{\mathcal{G} \times G,\mathcal{V} \oplus \mathbf{N}} \cong \mathcal{A}_{\mathcal{G},\mathcal{V} \oplus p^*\mathbf{N}}$, so the left-hand side of \eqref{equation-square} gives $p_*\mathcal{A}_{\mathcal{G},\mathcal{V} \oplus p^*\mathbf{N}}$.
    
    On the other hand, by \eqref{equation-fact1}, we have $\mathcal{A}_{\mathcal{G} \times G,\mathcal{V} \oplus \mathbf{N}} = \mathcal{A}_{\mathcal{G},\mathcal{V}} \boxtimes \mathcal{A}_{G,\mathbf{N}}$, and so the pushforward along the right map gives $\mathcal{R}_{G} \boxtimes \mathcal{A}_{G,\mathbf{N}}$ by \Cref{lemma-gt}. It follows that
    \begin{equation}\label{equation-final1}
        p_*\mathcal{A}_{\mathcal{G},\mathcal{V} \oplus p^*\mathbf{N}} \cong \mathcal{R}_{G} \otimes^! \mathcal{A}_{G,\mathbf{N}}.
    \end{equation}
    Applying the functor $H_{G(O)}^\bullet(-)$, the left-hand side yields the coordinate ring of $\mathcal{M}_C(\mathcal{G}, \mathcal{V} \oplus p^*\mathbf{N})$, while the right-hand side yields the coordinate ring of the relative Langlands dual $(T^*\mathbf{N})^\vee$. The same argument works with loop-rotation equivariance, by \cite[Proposition 3.3]{GannonWilliams2026}, and so the isomorphism is compatible with the Poisson structures. This finishes the proof of \Cref{thm-1}.

    \subsection{Proof of \texorpdfstring{\Cref{thm-2}}{Theorem 1.2}}

    We keep the notation from the introduction. Consider the diagram
    \[
        \begin{tikzcd}
            \mathcal{H} \times \mathcal{K} \arrow[r, "\iota_\phi"]\arrow[d, "q"'] & \mathcal{H} \times \mathcal{G} \arrow[d, "q \times p"] \\
            H \arrow[r, "\mathrm{graph}_{\phi}"] & H \times G
        \end{tikzcd}
    \]
    where $\mathrm{graph}_{\phi}$ is the graph embedding $h \mapsto (h,\phi(h))$, $\iota_\phi$ is as defined in \eqref{equation-homomorphism1}, and $q$ is short-hand for $q \circ \mathrm{pr}_{\mathcal{H}}$.

    We apply \eqref{equation-square} to the square above and the ring object $\mathcal{A}_{\mathcal{H} \times \mathcal{G}, \mathcal{U} \oplus \mathcal{V}}$. By \eqref{equation-fact2}, we have
    \[
        \iota_\phi^\dagger\mathcal{A}_{\mathcal{H} \times \mathcal{G}, \mathcal{U} \oplus \mathcal{V}} \cong \mathcal{A}_{\mathcal{H} \times \mathcal{K}, \iota_\phi^*(\mathcal{U} \oplus \mathcal{V})},
    \]
    so the left-hand side of \eqref{equation-square} gives $q_*\mathcal{A}_{\mathcal{H} \times \mathcal{K}, \iota_\phi^*(\mathcal{U} \oplus \mathcal{V})}$.

    Since $\mathcal{A}_{\mathcal{H} \times \mathcal{G}, \mathcal{U} \oplus \mathcal{V}} = \mathcal{A}_{\mathcal{H}, \mathcal{U}} \boxtimes\mathcal{A}_{\mathcal{G}, \mathcal{V}}$, it follows from \Cref{lemma-gt} that 
    \[
        (q \times p)_*\mathcal{A}_{\mathcal{H} \times \mathcal{G}, \mathcal{U} \oplus \mathcal{V}} \cong \mathcal{R}_{H \times G}.
    \]
    Thus \eqref{equation-square} in this case is
    \begin{equation}\label{equation-final2}
        q_*\mathcal{A}_{\mathcal{H} \times \mathcal{K}, \iota_\phi^*(\mathcal{U} \oplus \mathcal{V})} \cong (\mathrm{graph}_{\phi})^\dagger\mathcal{R}_{H \times G}.
    \end{equation}
    Taking $H(O)$-equivariant cohomology, the left-hand side yields the coordinate ring of the Coulomb branch $\mathcal{M}_C(\mathcal{H} \times \mathcal{K}, \iota_\phi^*(\mathcal{U} \oplus \mathcal{V}))$. By \cite[Theorem 2.2.2(ii)]{ginzburg2025pointwisepurityderivedsatake}, the $H(O)$-equivariant cohomology of the right-hand side is precisely the coordinate ring of $\overbar{\mathrm{Wh}_{H \to G}}$. This completes the proof of \Cref{thm-2}.
    
    \printbibliography

\end{document}